\documentclass[12pt,reqno]{amsart}
\allowdisplaybreaks[1]

\DeclareMathAlphabet{\mathpzc}{OT1}{pzc}{m}{it}

\usepackage{comment}
\usepackage{color,cancel,todonotes}
\usepackage[polish,english]{babel}
\usepackage{amsmath,amsfonts,amsxtra,fullpage}
\usepackage[normalem]{ulem}

\author{Hung Yean Loke}
\address{National University of Singapore, Science Drive 2, Singapore 117543}
\email{matlhy@nus.edu.sg}
\author{Tomasz Przebinda}
\address{Department of Mathematics, University of Oklahoma, Norman, OK 73019, USA}
\email{tprzebinda@ou.edu}

\title[The character and wave front set correspondence]{The character correspondence in the stable range \\ 
over a p-adic field}

\def\Chc{\mathrm{Chc}}

\def\g{\mathfrak g}

\def\sp{\mathfrak {sp}}

\def\cL{\mathcal{L}}

\def\Df{\Bbb{D}}

\def\fo{\mathfrak o}

\def\F{\mathbb{F}}
\def\R{\mathbb{R}}
\def\C{\mathbb{C}}

\def\ss1{\mathfrak s_{\overline 1}}

\def\hs1{\mathfrak h_{\overline 1}}

\def\supp{\mathrm{supp}}

\def\Op{\mathrm{Op}}

\def\G{\mathrm{G}}

\def\Zg{\mathrm{Z}}
\def\E{\mathrm{E}}
\def\Qg{\mathrm{Q}}

\def\H{\mathrm{H}}
\def\M{\mathrm{M}}

\def\L{\mathrm{L}}
\def\Bbb{\mathbb}

\def\A{\mathrm{A}}
\def\H{\mathrm{H}}
\def\T{\mathrm{T}}
\def\GL{\mathrm{GL}}

\def\SO{\mathrm{SO}}

\def\Sp{\mathrm{Sp}}

\def\Og{\mathrm{O}}

\def\fs{\mathfrak s}

\def \t{\tilde}
\def \wt{\widetilde}
\newcommand{\reg}[1]{ {#1}^{reg}}

\def\V{\mathsf{V}}
\def\Zv{\mathsf{Z}}
\def\W{\mathsf{W}}

\def\Vv{\mathrm{V}}
\def\Uv{\mathrm{U}}
\def\V{\mathsf{V}}

\def\X{\mathsf{X}}
\def\Y{\mathsf{Y}}
\def\Xv{\mathrm{X}}
\def\Yv{\mathrm{Y}}
\def\Ker{\mathrm{Ker}}

\def\End{\mathop{\hbox{\rm End}}\nolimits}

\def\det{\mathop{\hbox{\rm det}}\nolimits}

\def\Ad{\mathop{\hbox{\rm Ad}}\nolimits}
\def\Hom{\mathop{\hbox{\rm Hom}}\nolimits}

\def\tr{\mathop{\hbox{\rm tr}}\nolimits}

\def\lim{\mathop{\hbox{\rm lim}}\nolimits}

\def\supp{\mathop{\hbox{\rm supp}}\nolimits}

\def\Ss{\mathcal{S}}

\def\fontindex{\arabic}

\def\fonttitre{\textsf}
\newcounter{thh}

\newtheorem{thm}[thh]{\fonttitre{Theorem}}

\newtheorem{pro}[thh]{\fonttitre{Proposition}}
\newtheorem*{pro*}{\fonttitre{Proposition}}
\newtheorem{cor}[thh]{\fonttitre{Corollary}}
\newtheorem*{coro*}{\fonttitre{Corollary}}
\newtheorem{lem}[thh]{\fonttitre{Lemma}}

\newtheorem*{defi*}{\fonttitre{Definition}}
\newtheorem{conj}{\fonttitre{Conjecture}}

\newtheorem*{nota*}{\fonttitre{Notation}}
\newenvironment{prf}{\begin{proof}}{\end{proof}}
\def\muet{ \ifthenelse{\equal{a}{b}}}
\def\nn{\nonumber}
\def\Z'{\Bbb{Z}'}

\def\biblio{\sloppy 
\bibliographystyle{alpha}            
\bibliography{article}}

\begin{document}
\thanks{The first author is grateful to the University of Oklahoma for hospitality and financial support in February 2017. The second author gratefully acknowledges hospitality and financial support from the Institute of Mathematical Sciences at the National University of Singapore and the National Science Foundation under Grant DMS-2225892. }

\date{}



\maketitle

\begin{abstract}
Given a real irreducible dual pair there is  an integral kernel operator 
which maps the distribution character  of an irreducible admissible representation  of the group with the smaller or equal rank to an invariant 
eigen\-distribution  on the group with the larger or equal rank. If the pair is in the stable range  and if the representation is unitary, then the resulting distribution is the character of the representation obtained via Howe's correspondence. This construction was transferred to the p-adic case and a conjecture was formulated.

In this note we verify a weaker version of this conjecture for dual pairs in the stable range over a p-adic field. 
\end{abstract}

\tableofcontents

\section{\bf Introduction. \rm}\label{Introduction}
For a real irreducible dual pair $(\G,\G')$ with the rank of $\G'$ less or equal to the rank of $\G$ \cite{PrzebindaCauchy} provides an integral kernel operator $\Chc$
 which maps the distribution character $\Theta_{\Pi'}$ of an irreducible admissible representation $\Pi'$ of $\wt\G'$ to an invariant 
eigendistribution $\Theta'_{\Pi'}$ on the group $\wt\G$ with the correct infinitesimal character, \cite{BerPrzeCHC_inv_eig}. If the pair is in the stable range with $\G'$ the smaller member and if the representation is unitary, then $\Theta'_{\Pi'}=\Theta_\Pi$, where $\Pi$ is associated to $\Pi'$ via Howe's correspondence, \cite{PrzebindaStableUnitary}.  The acronym $\Chc$ stands for the Cauchy Harish-Chandra integral, because as explained in \cite{PrzebindaCauchy} the construction gives a direct link from the Cauchy determinantal identity through Harish-Chandra's theory of the semisimple orbital integrals to Howe's correspondence. This construction was transferred to the $p$-adic case in \cite{LokePrzebinda_chc_padic_def}.

In this note we verify a weaker version of the conjecture formulated in \cite{LokePrzebinda_chc_padic_def} for dual pairs in the stable range over a $p$-adic field. 
Let $Z'\subseteq \G'$ denote the center and let  $\G'{}^\circ\subseteq \G'$ be the Zariski identity component of $\G'$. Then $Z'\G'{}^\circ=\G'$ unless $\G'$ is an even orthogonal group.

\begin{thm}\label{theteequalstheta'}
Suppose $(\G, \G')$ is an irreducible dual of type I in the stable range with $\G'$ - the smaller member. 
Then there is a non-empty Zariski open subset $\G''\subseteq \G$ with the following property.  

Let $\Pi'$ be any genuine irreducible unitary representation of $\wt\G'$ and let $\Pi$ be the representation of $\wt\G$ corresponding to $\Pi'$. Let $\Theta_\Pi$ denote the distribution character of $\Pi$. Recall the distribution $\Theta_{\Pi'}'$ on $\wt \G$, \cite[(130)]{LokePrzebinda_chc_padic_def}. 
Assume that the character
$\Theta_{\Pi'}$ of the representation $\Pi'$ is supported in $Z'\G'{}^\circ$.
Then 
\[
\Theta_\Pi(\Psi)=\Theta_{\Pi'}'(\Psi) \qquad (\Psi\in C_c(\wt \G''))\,.
\]
\end{thm}
The proof follows the argument used in \cite{PrzebindaStableUnitary}.

\section{\bf The  Weil representation. \rm}\label{The  Weil representation}
In this section we recall the Weil representation \cite{WeilWeil} with the details suitable for our computations following \cite{AubertPrzebinda_omega}.
 Fix a non-trivial unitary character 
$\chi \colon \F\to\C^\times$ of the additive group $\F$ with the kernel equal to $\fo_\F$, the ring of integers in $\F$. We assume that the Haar measure of $\fo_\F$ is $1$. 

Let $\W$ be a finite dimensional vector space over $\F$ with a non-degenerate symplectic form $\langle \cdot,\cdot \rangle$. 
Fix a lattice $\cL\subseteq \W$ and the corresponding norm 
\[
N_\cL(w)=\inf\{|a|^{-1}:\ a\in \F^\times,\ aw\in \cL\} \qquad (w\in \W)\,.
\]
We shall assume that the lattice $\cL$ is self-dual in the sense that
\[
\left(\langle w,w'\rangle\in\fo_\F\ \ \text{for all}\ \ w'\in \mathcal L\right)\  \Longleftrightarrow\ w\in \mathcal L\,.
\]
For any subspace $\Uv\subseteq \W$ we normalize the Haar measure $\mu_\Uv$ on $\Uv$ so that the volume of the lattice $\cL\cap \Uv$ is $1$.
If $\Vv\subseteq \Uv$, then  we  normalized Haar measure $\mu_{\Uv/\Vv}$  so that the volume of the lattice $(\cL\cap\Uv+\Vv)/\Vv$ is $1$.
If 
$q$ is a nondegenerate quadratic form on $\Uv$, then we set
\[
\gamma(q)=\underset{r\to \infty}{\lim}\,\int_{u\in \Uv, |u|<r}\chi(\frac{1}{2}q(u))\,d\mu_\Uv(u)\ \ \ 
\text{and}\ \ \ 
\gamma_{\mathrm {Weil}}(q)=\frac{\gamma(q)}{|\gamma(q)|}\,.
\]
Here $\chi(\frac{1}{2}q(u))$ is a Gaussian and $\gamma_{\mathrm {Weil}}(q)$ is the Weil factor of $q$, with $\gamma_{\mathrm {Weil}}(q)^8=1$. 
In particular if $\Uv=\F$ and $\cL \cap U = \fo_F$ we have
\[
\gamma(a)=\underset{r\to \infty}{\lim}\,\int_{u\in \F, |u|<r}\chi(\frac{1}{2}au^2)\,du\ \ \ 
\text{and}\ \ \ 
\gamma_{\mathrm {Weil}}(a)=\frac{\gamma(a)}{|\gamma(a)|}\qquad (a\in\F^\times)\,.
\]
The symplectic group $\Sp(\W)\subseteq \End(\W)$ is the group of the isometries of the symplectic form~$\langle \cdot,\cdot \rangle$. 
Define
\[
\det\left(g-1 \colon \W/\Ker\,(g-1) \to (g-1)\W \right) (\fo_\F^\times)^2
=\det(\langle (g-1) w_i, w_j\rangle_{1\leq i,j\leq m})(\fo_\F^\times)^2\in \F^\times/(\fo_\F^\times)^2\,,
\]
where $w_1$, $\dots$, $w_m$ are such that
\[
\W = \F w_1 \oplus \dots \oplus \F w_m \oplus \Ker(g-1)
\]
and the summands on the right are $N_\cL$-orthogonal, i.e.
\[
N_\cL(a_1w_1+\cdots+a_mw_m+w)=\max\{N_\cL(a_1w_1), \dots, N_\cL(a_mw_m), N_\cL(w)\}\,.
\]
For $g, g_1, g_2\in \Sp$,  let
\[
\Theta^2(g)=\gamma(1)^{2\dim\,(g-1)\W-2}\,\big[\gamma(\det (g-1 \colon \W/\Ker(g-1)\to (g-1)\W))\big]^2
\]
\[
C(g_1,g_2)=\sqrt{\left|\frac{\Theta^2(g_1g_2)}{\Theta^2(g_1)\Theta^2(g_2)}\right|}\;
\gamma_{\mathrm{Weil}}(q_{g_1,g_2}),
\]
where
\vskip-9mm
\begin{multline*}
q_{g_1,g_2}(u',u'')=\frac{1}{2}\langle (g_1+1)(g_1-1)^{-1}u',u''\rangle\\
+\frac{1}{2}\langle (g_2+1)(g_2-1)^{-1}u',u''\rangle\\ 
(u',u''\in (g_1-1)\W\cap (g_2-1)\W)\,.
\end{multline*}
The Metaplectic group is defined as
\[
\wt{\Sp}=\left\{\t g=(g,\xi)\in \Sp\times \C,\ \ \xi^2=\Theta^2(g)\right\}
\]
with the multiplication
\[
(g_1,\xi_1)(g_2,\xi_2)=(g_1g_2,\xi_1\xi_2 C(g_1,g_2))\,.
\]
Let $\W=\Xv\oplus \Yv$ be a complete polarization. Set
\[
\Op \colon \Ss^*(\Xv\times \Xv)\to \Hom_{ \C}(\Ss(\Xv),\Ss^*(\Xv))\,,\ \ \ 
\Op(K)v(x)=\int_\Xv K(x,x')v(x')\,d\mu_\X(x').
\]
Recall the Weyl transform $\mathcal K \colon \Ss^*(\W)\to \Ss^*(\Xv\times \Xv)$,
\[
\mathcal K(f)(x,x')=\int_\Yv f(x-x'+y)\chi\big(\tfrac{1}{2}\langle y, x+x'\rangle\big)\,d\mu_\Y(y)
\]
and an imaginary Gaussian on $(g-1)\W$
\[
\chi_{c(g)}(u)=\chi\big(\tfrac{1}{4}\langle \underbrace{(g+1)(g-1)^{-1}}_{\scriptsize c(g)}u, u\rangle\big) \qquad (u=(g-1)w,\ w\in\W)\,.
\]
Let
\begin{equation} \label{eqrho}
\rho=\Op\circ \mathcal K \colon \Ss^*(\W)\to \Hom_{\C}(\Ss(\Xv),\Ss^*(\Xv)) \,.
\end{equation}
For $\t g=(g,\xi)\in\wt\Sp$, we define
\[
\Theta(\t g)=\xi,\qquad T(\t g)=\Theta(\t g)\chi_{c(g)}\mu_{(g-1)\W},\qquad 
\omega(\t g)=\rho \circ T(\t g)\,.
\]
One could deduce from Lemma 5.11 in \cite{AubertPrzebinda_omega} that $(\omega, \Ss(\Xv))$ is a representation of $\wt\Sp(\W)$. In addition by Theorem 5.26 in \cite{AubertPrzebinda_omega} $(\omega, \L^2(\Xv))$ is the Schr\"odinger model of the Weil representation of $\wt\Sp(\W)$ attached to the character $\chi$.  Furthermore, the cocycle
\[
C(g_1,g_2)=\frac{\Theta(\t g_1\t g_2)}{\Theta(\t g_1)\Theta(\t g_2)} \qquad (\t g_1, \t g_2\in \wt{\Sp}(\W))\,.
\]
It turns out that $\Theta$ is the distribution character of $\omega$. We shall refer to $T(\t g)$  as a normalized Gaussian. For future reference we notice that
\begin{equation} \label{eqtromegarho}
\tr \omega(\t g)\rho(\phi)=T(\t g)(\phi) \qquad (\phi\in \Ss(\W))\,.
\end{equation}
Indeed, in terms of generalized functions, the left hand side is equal to
\[
\tr \rho(T(\t g))\rho(\phi)=\int_\X\int_\X \mathcal K(T(\t g))(x,x')\mathcal K(\phi)({x',x})\,dx\,dx'
\]
which is equal to the right hand side because of the definition of $\mathcal K$ and the Fourier inversion formula. 

\section{\bf A mixed model  of the  Weil representation. \rm}\label{A mixed model of the  Weil representation}
 For a subset $S$ of $\Sp(\W)$, we denote its inverse image in $\wt{\Sp}(\W)$ by $\t S$. For $g \in \Sp(\W)$, we denote an element of the inverse image of $g$ in $\wt \Sp(\W)$ by $\wt g$.

In this section we recall the explicit  formulas for $\omega(\t g)$ for some particular elements 
$\t g\in \wt{\Sp}(\W)$. Recall the function
\[ 
\F^\times\ni a\to \fs(a):= |a|_\F \,\frac{\gamma(a)^2}{\gamma(1)^2}=\frac{\gamma_{\mathrm{Weil}}(a)^2}{\gamma_{\mathrm{Weil}}(1)^2}\in\C^\times\,.
\] 
This is a unitary character of $\F^\times/(\F^\times)^2$.

For a subset $\M\subseteq \End(\W)$ let $\M^c=\{m\in \M\,:\,\det (m-1)\ne 0\}$ denote the domain of the Cayley transform in $\M$.
\begin{pro}\label{formula for M R}
Let $\M\subseteq \Sp(\W)$ be the subgroup of all the elements that preserve 
$\X$ and $\Y$. Set
\[
\det_\X^{-1/2}(\t m)=\Theta(\t m) |\det(\frac{1}{2}(c(m|_\X)+1))|^{-1} \qquad (\t m\in \wt\M^c).
\]
Then
\[
\left(\det_\X^{-1/2}(\t m)\right)^2=\fs(\det(m|_\X))^{-1}|\det(m|_\X)|^{-1} \qquad (\t m\in \wt\M^c)\,,
\]
the function $\det_\X^{-1/2}\colon
\wt\M^c\to \C^\times$ extends to a continuous group homomorphism
\[
\det_\X^{-1/2}\colon\wt\M\to \C^\times.
\]
For $\t m\in\wt\M$ and $v\in \Ss(\X)$, we have $\omega(\t m)v \in \Ss(\X)$. It is given by
\[
\omega(\t m)v(x)=\det_\X^{-1/2}(\t m) v(m^{-1}x) \qquad (x\in\X).
\]
\end{pro}

Suppose $\W=\W_1\oplus\W_2$ is the direct orthogonal sum of two symplectic spaces. There are inclusions 
\begin{equation}\label{embedingofproductofgroups0}
\Sp(\W_1)\subseteq \Sp(\W),\ \ \ \Sp(\W_2)\subseteq \Sp(\W)
\end{equation}
defined by
\begin{eqnarray*}
g_1(w_1+w_2)&=&g_1w_1+w_2\,\\
g_2(w_1+w_2)&=&w_1+g_2w_2  \qquad (g_j\in \Sp(\W_j),\ w_j\in \W_j,\ j=1,2)\,.
\end{eqnarray*}
Furthermore, the map
\begin{equation}\label{embedingofproductofgroups}
\Sp(\W_1)\times \Sp(\W_2)\ni (g_1, g_2)\to g_1g_2\in \Sp(\W)
\end{equation}
is an injective group homomorphism.

Assume that $\cL\cap\W_1\oplus\cL\cap\W_2=\cL$. Then we have two metaplectic groups $\wt\Sp(\W_j)$, $j=1,2$. (Here $\wt\Sp(\W_j)$ is defined using the same $\Theta^2$.) It is not difficult to see that the embeddings \eqref{embedingofproductofgroups0} lift to the embeddings
\[
\wt\Sp(\W_1)\subseteq \wt\Sp(\W),\ \ \ \wt\Sp(\W_2)\subseteq \wt\Sp(\W).
\]
It follows easily from the formula for the cocycle that
\[
C(g_1, g_2)=1 \qquad (g_j\in \Sp(\W_j),\ j=1,2)\,.
\]
Hence \eqref{embedingofproductofgroups} lifts to a group homomorphism
\[
\wt\Sp(\W_1)\times \wt\Sp(\W_2)\ni (\t g_1, \t g_2)\to \t g_1\t g_2\in \wt\Sp(\W)\,,
\]
with kernel equal to a two-element group. Moreover, in terms of the identification
\[
\Ss(\W)=\Ss(\W_1)\otimes \Ss(\W_2)\,,
\]
we have
\[
T(\t g_1\t g_2)=T_{1}(\t g_1)\otimes T_{2}(\t g_2) \qquad (\t g_j\in \wt\Sp(\W_j),\ j=1,2)\,,
\]
where $T_{j}(\t g_1)$ is the normalized Gaussian for the space $\W_j$,  $j=1,2$. Hence,
\[
\omega(\t g_1\t g_2)=\omega_{1}(\t g_1)\otimes \omega_{2}(\t g_2) \qquad (\t g_j\in \wt\Sp(\W_j),\ j=1,2)\,,
\]
where $\omega_j$ is the Weil representation of $\wt\Sp(\W_j)$ for $j=1,2$.

Suppose from now on that $\W_j=\X_j\oplus \Y_j$, $j=1,2$, are complete polarizations such that
\[
\X=\X_1\oplus\X_2\ \ \ \text{and}\ \ \ \Y=\Y_1\oplus \Y_2.
\]
Then, in particular, we have the following identifications
\begin{equation}\label{tensorproductidentification}
\Ss(\X)= \Ss(\X_1)\otimes  \Ss(\X_2)= \Ss(\X_1,\Ss(\X_2)).
\end{equation}
\begin{cor}\label{Maction}
Let $m_j \in \Sp(W_j)$ for $j = 1,2$. We assume that $m_1$ preserves $\X_1$ and $\Y_1$.
Then for $v_1\in \Ss(\X_1)$, $v_2\in \Ss(\X_2)$, $x_1\in \X_1$ and $x_2\in \X_2$,
\[
\left(\omega(\wt {m_1}\wt {m_2})(v_1\otimes v_2)\right)(x_1+x_2)=\det_{X_1}^{-1/2}(\wt{m_1})v_1(m_1^{-1}x_1)(\omega_2(\wt{m_2}) v_2)(x_2)\,.
\]
Thus, in terms of \eqref{tensorproductidentification},
\[
\omega(\wt {m_1}\wt {m_2})v(x_1)=\det_{X_1}^{-1/2}(\wt{m_1})\omega_2(\wt{m_2})v(m_1^{-1}x_1) \qquad (v\in \Ss(\X_1,\Ss(\X_2)),\ x_1\in \X_1)\,.
\]
\end{cor}
\begin{pro}\label{N1action}
Suppose $n\in \Sp(\W)$ acts trivially on $\Y_1^\perp(=\Y_1+\W_2)$. 
Then for
$v\in \Ss(\X_1,\Ss(\X_2))$ and $x_1\in \X_1$,
\[
\omega(\t n)v(x_1)=\pm \chi_{c(-n)}(2x_1)v(x_1)\,.
\]
\end{pro}
Proposition \ref{N1action} is well known. If $\W_2=0$ then it coincides with \cite[Proposition 5.28]{AubertPrzebinda_omega}. The general case may be verified via an argument similar to the one used there.

\section{\bf The Cauchy Harish-Chandra integral.}


In this section we recall some results and a conjecture in \cite{LokePrzebinda_chc_padic_def}.

For any $\Psi\in C_c^\infty(\wt{\Sp}(\W))$ the distribution
\[
T(\Psi)=\int_{\wt{\Sp}(\W)} \Psi(g) T(g) \,dg\in \Ss^*(\W)
\]
is a function that belongs to $\Ss(\W)$ (times the measure $dw$) and the formula
\[
\Chc(\Psi)=\int_\W T(\Psi)(w)\,dw
\]
defines a distribution of $\wt{\Sp}(\W)$. 

\medskip

Let $\Bbb D$ be a division $\F$-algebra, with a possibly trivial involution $\sigma$ fixing $\F$ pointwise. 
From now on, all $\Bbb D$-modules are right $\Bbb D$-modules unless otherwise stated. For two $\Bbb D$-modules $V_1$ and $V_2$, $\Hom_{\Bbb D}(V_1,V_2)$ denotes the set of right $\Bbb D$-module homomorphisms.

Let $\epsilon = \pm 1$. Let $\V$ and $\V'$ be free $\Bbb D$-modules of finite rank. Let $(\cdot,\cdot)$ be a non-degenerate $\epsilon$-hermitian form on $\V$ and let $(\cdot,\cdot)'$ be a nondegenerate $(-\epsilon)$-hermitian form on $\V'$. We set
\begin{eqnarray}
\G & = & \{ g \in \End_{\Bbb D}(\V) : (gu,gv) = (u,v) \text{ for all } u, v \in \V \} \ \ \text{ and } \label{eqG} \\
\G' & = & \{ g' \in \End_{\Bbb D}(\V') :  (g'u',g'v')' = (u',v')' \text{ for all } u', v' \in \V' \}. \nn
\end{eqnarray}
Then $\G \cdot \G'\subseteq \Sp(\W)$ is an irreducible dual pair. 
Let $\H'\subseteq \G'$ be a Cartan subgroup with the split part $\A'\subseteq \H'$. Let $\A''$ denote the centralizer of $\A'$ in $\Sp(W)$. 
It is shown in \cite{LokePrzebinda_chc_padic_def} that $(\A'',\A')$ is a dual pair.
We prove the following lemma in \cite[Section 4]{LokePrzebinda_chc_padic_def}.

\begin{lem}\label{question2Chc}
For any $\Psi\in C_c^\infty(\wt{\A''{}^c})$, the distribution
\begin{equation}\label{theChc0GLChc}
T(\Psi)=\int_{\wt{\A''{}^c}} \Psi(\t g)T(\t g)\,d\t g\in \Ss^*(\W)
\end{equation}
is a function on $\W$. The formula
\begin{equation}\label{theChc2GLChc}
\Chc(\Psi)=\int_{\A'\backslash\W_{\A'}} T(\Psi)(w)\,d(\A' w)
\qquad (\Psi\in C_c^\infty(\wt{\A''{}^c}))
\end{equation}
defines a distribution on $\wt{\A''{}^c}$ which coincides with a complex valued measure.
This measure extends by zero to $\wt{\A''}$ and defines a distribution $\Chc$ on $\t \A''$.
\end{lem}

Let $\t h'\in \reg{\wt{\H'}}$. We define an embedding $\iota \colon \wt \G \rightarrow \wt \A''$ by $\iota(g) = gh'$. 
The pullback of the distribution $\Chc$ via $\iota$ to $\wt{\G}$ is well-defined. 
We will denote this pullback distribution on $\wt \G$ by $\Chc_{\t h'}$.

Let $\Pi'$ be an irreducible admissible representation of $\wt{\G'}$ which occurs in Howe's correspondence for the pair $(\wt\G, \wt{\G'})$ and let $\Pi_1$ be the corresponding maximal Howe's quotient  representation of $\wt\G$. Let $\Theta_{\Pi'}$ denote the distribution character of $\Pi'$. We state the following conjecture in 
\cite[Section 5]{LokePrzebinda_chc_padic_def}.
\begin{conj}\label{Conjectureforcharacters}
For $\Psi\in C_c(\wt\G)$, we set
\begin{equation}\label{Conjectureforcharacters2}
\Theta_{\Pi'}'(\Psi)=
C_{\Pi'}\sum_{\H'}\frac{1}{|W(\H')|}\int_{\reg{\H'}}\check\Theta_{\Pi'}(\t h')
|\Delta(h')|^2\frac{1}{\mathrm{volume}(\A'\backslash\H')}
\Chc_{\t h'}(\Psi)\,d\t h'\,,
\end{equation}
where $\check\Theta_{\Pi'}(\t h')=\Theta_{\Pi'}(\t h'^{-1})$,  $\Delta(h')$ is the Weyl denominator and
\[
C_{\Pi'}=(\text{the central character of $\Pi'$ evaluated at $\wt{-1}$})^{-1}\cdot\Theta(\wt{-1})\,.
\]
Then $\Theta_{\Pi'}'$ is a distribution on $\t G$. Moreover 
\begin{equation}\label{Conjectureforcharacters3}
\Theta_{\Pi'}'=\Theta_{\Pi_1}\,.
\end{equation}
as distributions.
\end{conj}

Our goal in this paper is to prove the above conjecture when $(\G,\G')$ is in stable range where $\G'$ is the smaller member.

\section{\bf The $p$-adic method of stationary phase. \rm}\label{stationary phase}
In this section we recall \cite[Proposition 1.1]{Heifetz}, or rather more numerically explicit \cite[Lemma A.14]{LokePrzebinda_chc_padic_def}, in a coordinate free formulation. Let $\Zv$ be a finite dimensional vector space over $\F$ with a norm $|\cdot|$ and let $\Zv^*$ be the dual vector space with the corresponding norm denoted by the same symbol 
\[
|z^*|=\max_{z\in\Zv\,,\ |z|=1}|z^*(z)| \qquad (z^*\in\Zv^*)\,.
\]
\begin{lem}\label{Theorem 7.7.1Hormander}\cite[Theorem 7.7.1]{Hormander}, \cite[Proposition 1.1]{Heifetz}
Let $U\subseteq \Zv$ be an open compact subset and let $f:U\to\F$ be a differentiable function such that
\begin{equation}\label{Theorem 7.7.1Hormander1}
f(z_0+z)=f(z_0)+f'(z_0)(z)+R(z_0,z)(z)(z) \qquad (z_0,\ z_0+z\in U)\,,
\end{equation}
where $f'(x_0)\in\Zv^*$ is the derivative of $f$ at $z_0$ and $R(z_0,z)\in\Hom(\Zv,\Zv^*)$ is a linear function with 
\begin{equation}\label{Theorem 7.7.1Hormander2}
M:=\max_{z_0, z_0+z\in U, |y|=1}|R(z_0, z)(y)(y)|<\infty\,.
\end{equation}
We also assume that
\begin{equation}\label{Theorem 7.7.1Hormander3}
\delta:=\min_{z_0\in U}|f'(z_0)|>0\,.
\end{equation}
Denote by $B_r\subseteq \Zv$ the closed ball centered at $0$, with radius $r$. 
Let $\phi\in \Ss(U)$ and
let $m_0\in \Bbb Z$ be the minimum of the $m\in \Bbb Z$ such that there is a finite disjointed covering
\begin{equation}\label{Theorem 7.7.1Hormander4}
U=\bigsqcup_k \left(z_k+B_{q^{-m}}\right)
\end{equation}
and $\phi$ is constant on each $z_k+B_{q^{-m}}$.
(The covering exists because $\phi$ is locally constant and $U$ is open and compact.)
Then
\begin{equation}\label{Theorem 7.7.1Hormander5}
\int_U\chi(\lambda f(z))\phi(z)\,dx=0 \qquad(\lambda \in \F^\times\,,\ |\lambda|>\max\{q\delta^{-2}M, \delta^{-1}q^{m_0}\})\,.
\end{equation}
\end{lem}
\section{\bf The restriction of the Weil representation to the dual pair. \rm}\label{The restriction of the Weil representation to the dual pair}
Let $\Bbb D$ be a division $\F$-algebra, with a possibly trivial involution $\sigma$ fixing $\F$ pointwise. We recall the right $\Bbb D$-modules  $\V$, $\V'$ and the irreducible dual pair $(\G, \G')$ in \eqref{eqG}. We suppose $(\G, \G')$ is in stable range where $\G'$ is the smaller member.

The stable range assumption means that there is an isotropic subspace $\X_{(1)}\subseteq \V$ such that $\dim \V'\leq \dim \X_{(1)}$. We select an isotropic subspace $\Y_{(1)}\subseteq \V$,  complementary to $\X_{(1)}^\perp$. Let $\V_{(2)}\subseteq \V$ be the orthogonal complement of $\X_{(1)}\oplus \Y_{(1)}$, so that $\V=\X_{(1)}\oplus\V_{(2)}\oplus Y_{(1)}$.

The symplectic space  may be realized as $\W=\Hom_{\Bbb D}(\V,\V')$ with 
\begin{equation}\label{explicitsymplecticform}
\langle w',w\rangle=\tr_{\Df/\F}(w^*w'),
\end{equation}
where $w^*\in \Hom_{\Bbb D}(\V',\V)$ is defined by
$(wv,v')'=(v,w^*v')$, with $v\in \V$ and  $v'\in \V'$. The group $\G'$ acts on $\W$ by the post-multiplication and the group $\G$ by the pre-multiplication by the inverse.
Set 
\[
\X_1=\Hom_{\Bbb D}(\X_{(1)},\V'), \ \ \Y_1=\Hom_{\Bbb D}(\Y_{(1)},\V') \ \text{ and } \ \W_2=\Hom_{\Bbb D}(\V_{(2)},\V').
\] 
Then $\Y_1$ and $\X_1^\perp$ are complementary subspaces of $\W$.  With respect to the symplectic form~\eqref{explicitsymplecticform}, $\W_2$ is the orthogonal complement of $\W_1=\X_1+\Y_1$. 
We shall work in the mixed model of the Weil representation adapted to the decomposition $\W=\X_1\oplus \W_2\oplus \Y_1$, as explained in the  Section \ref{A mixed model of the  Weil representation}.

Denote by $i_{\Y_1}:\Y_1\to\X_1\oplus\W_2\oplus \Y_1$ the injection and by  $p_{\X_1}:\X_1\oplus\W_2\oplus \Y_1\to\X_1$ the projection. 
In particular, for  $z\in\sp(\W)$, we have a linear map $p_{\X_1}zi_{\Y_1} \colon \Y_1\to\X_1$. If $z\in\g$ then the above map is bijective if and only if
the map  $p_{\X_{(1)}}zi_{\Y_{(1)}} \colon \Y_{(1)}\to\X_{(1)}$ is bijective, where $i_{\Y_{(1)}}:\Y_{(1)}\to\X_{(1)}\oplus\V_2\oplus \Y_{(1)}$ the injection and by  $p_{\X_{(1)}}:\X_{(1)}\oplus\V_2\oplus \Y_{(1)}\to\X_{(1)}$ the projection. There is one case when there is no $z\in\g$, such that $p_{\X_{(1)}}zi_{\Y_{(1)}}$ is bijective. This happens if $\G$ is an orthogonal group (i.e. the involution $\sigma$ is trivial and the form $(\cdot,\cdot)$ is symmetric) and the dimension of $\X_{(1)}$ is odd. 
By the stable range assumption $\dim\X_{(1)}\geq \dim \V'$. Hence we may choose $\X_{(1)}$ with $\dim\X_{(1)}=\dim \V'$, which is even.
Thus the set of elements $z\in\g$ such that the map $p_{\X_1}zi_{\Y_1} \colon \Y_1\to\X_1$ is bijective, or equivalently the form 
$
q_{z,\Y_1}(y,y')=\frac{1}{2}\langle z y,y'\rangle
$ 
is non-degenerate on $\Y_1$, is not empty.  We shall fix such a choice, where $p_{\X_1}zi_{\Y_1} \colon \Y_1\to\X_1$ is bijective, for the rest of this article.

(By the way notice that up to this point, we need stable range if $\G$ is an orthogonal group and $\dim \X_{(1)}$ is odd. 
Then we change $\X_{(1)}$ such that $\dim\X_{(1)}=\dim \V'$ is even. 
For other groups $\G$, $\X_{(1)}$ is arbitrary and there is no need for stable range. Of course we need stable range for other reasons later.)

The complete polarization $\W_1=\X_1\oplus \Y_1$ leads to the Weyl transform $\mathcal K_1 \colon \Ss^*(\W_1)\to\Ss^*(\X_1\times\X_1)$. 
Hence $\mathcal K_1\otimes 1 \colon \Ss^*(\W)\to\Ss^*(\X_1\times\X_1\times \W_2)$. In order to shorten the notation we shall write $\mathcal K_1$ for $\mathcal K_1\otimes 1$. Explicitly
\[
{\mathcal{K}}_1(f)(x,x', w_2)=\int_{\Y_1} f(x-x'+y+w_2)\chi\big(\frac{1}{2}\langle y, x+x'\rangle\big)\,dy \qquad (f \in \Ss(\W), x,x'\in \X_1, w_2\in \W_2)\,.
\]
(This is a function on $\X_1 \times \X_1 \times \W_2$.)
For $z \in \End(\W)$, we define $\chi_z(w) = \chi(\frac{1}{4}\langle zw,w \rangle)$ for $w \in \W$.
\begin{lem}\label{partialweyltransformfor g}
Let $z\in \g^c$ be such that $p_{\X_1}zi_{\Y_1}$ is invertible, with the inverse 
\[
(p_{\X_1}zi_{\Y_1})^{-1} \colon \X_1\to\Y_1\,. 
\]
Then for $x,x'\in\X_1$ and $w_2\in \W_2$ we have
\begin{eqnarray}\label{partialweyltransformfor g1}
\mathcal K_1(T(\wt{c(z)}))(x,x',w_2)&=&\Theta(\wt{c(z)})\gamma(q_{z,\Y_1})\\
&&\chi_z(x-x')\chi_{(p_{\X_1}zi_{\Y_1})^{-1}}(x+x'-p_{\X_1}(z(x-x')+zw_2))\nn\\
&&\chi(\frac{1}{2}\langle zw_2,x-x'\rangle)\chi_z(w_2)\,.\nn
\end{eqnarray}
Let $h\in \G$ be the element that acts via multiplication by $-1$ on $\W_1$ and by the identity on $\W_2$. Suppose that in addition $\det(hc(z)-1)\ne 0$ and let $z_h=c(hc(z))$. Then
\[
\mathcal K_1(T(\wt{c(z_h)}))(x,x',w_2)=\det_{\X_1}^{-1/2}(\t h)\mathcal K_1(T(\wt{c(z)}))(x,-x',w_2)\,.
\]
Equivalently,
\[
\mathcal K_1(T(\wt{c(z)}))(x,x',w_2)=\det_{\X_1}^{1/2}(\t h)\mathcal K_1(T(\wt{c(z_h)}))(x,-x',w_2)\,.
\]
(Here $\t h$ is one of the two elements in the preimage of $h$  chosen so that the right hand side is equal to the left hand side.)
\end{lem}
\begin{prf}
This is verified by the argument used to prove Proposition 5.29 in \cite{AubertPrzebinda_omega} applied to an element $g_1\in \Sp(\W_1)$ such that 
$g_1$ acts trivially on $\X_1$ and $\W_1/\X_1$ and the restriction of $c(-g_1)$ to $\Y_1$ is equal to $p_{\X_1}zi_{\Y_1}$. 
\end{prf}
Here is a technical lemma, analogous to \cite[Lemma 4.3]{DaszPrzebindaInv}. Recall that for a test function $\Psi\in C_c^\infty(\wt\G)$
\[
T(\Psi)=\int_{\wt\G}\Psi(g)T(g)\,dg
\]
is a well defined distribution on $\W$. Hence $\mathcal K_1(T(\Psi))$ is a tempered distribution on $\X_1\times \X_1\times \W_2$. Fix a norm $|\cdot|$ on the $\F$ vector space $\End(\V)$, see \cite[Definition 1, chapter II, paragraph 1]{Weil_Basic}. 
We may assume that
\begin{equation}\label{operatornorminequality}
|z_1z_2|\leq |z_1||z_2| \qquad (z_1, z_2\in \End(\V))\,.
\end{equation}
Given $x \in \X_1 = \Hom_{\Bbb D} (\X_{(1)}, \V')$, we extend $x$ trivially over $\Y_{(1)} \oplus \V_{(2)}$ so $x \in \Hom_{\Bbb D}(\V,\V')$. In particular we have $x^* \in \Hom_{\Bbb D}(\V',\V)$. 
\begin{lem}\label{maintechicallemma}
There is a Zariski open subset $\G''\subseteq \G$ such that for $\Psi\in C_c^\infty(\wt\G'')$
the distribution $\mathcal K_1(T(\Psi))$ is a locally constant function on $\X_1\times \X_1\times \W_2$.
Moreover, there is a constants $C_\Psi$ such that 
\begin{eqnarray}\label{maintechicallemma1}
&&\mathcal K_1(T(\Psi))(x,x',w_2)=0\ \ \  \text{if}\\
&&|x^*x|+|x'{}^*x'|+|x^*x'|+|x'{}^*x|+|x^*w_2|+|x'{}^*w_2|+|w_2^*w_2|>C_\Psi\,.\nn
\end{eqnarray}
\end{lem}
\begin{prf}
The function \eqref{partialweyltransformfor g1} is equal to the  locally constant function $\Theta(\wt{c(z)})\gamma(q_{z,\Y_1})\ne 0$, times $\chi(\frac{1}{4}\phi_{x,x',w_2}(z))$, where
\begin{eqnarray*}
\phi_{x,x',w_2}(z)&=&\langle z(x-x'),x-x'\rangle\\
&+&\langle (p_{\X_1}z i_{\Y_1})^{-1}(x+x'-p_{\X_1}(z(x-x')+zw_2)), x+x'-p_{\X_1}(z(x-x')+zw_2)\rangle\\
&+&2\langle zw_2,x-x'\rangle+\langle zw_2,w_2\rangle\,.
\end{eqnarray*}
In order to simplify the computations we introduce the following notation
\begin{eqnarray*}
&& A(z)=p_{\X_{(1)}}zi_{\X_{(1)}},\ \ \ B(z)=p_{\X_{(1)}}zi_{\Y_{(1)}},\ \ \ C(z)=p_{\Y_{(1)}}zi_{\X_{(1)}},\ \ \ F(z)=C(z)^{-1},\\ 
&&D(z)=p_{\V_{(2)}}zi_{\Y_{(1)}},\ \ \ E(z)=p_{\V_{(2)}}zi_{\X_{(1)}},\ \ \ z_2=p_{\V_{(2)}}zi_{\V_{(2)}}\,,
\end{eqnarray*}
where
\[
i_{\X_{(1)}} \colon \X_{(1)}\to \X_{(1)}\oplus \V_{(2)}\oplus \Y_{(1)}\,,\ 
i_{\Y_{(1)}} \colon \X_{(1)}\to \X_{(1)}\oplus \V_{(2)}\oplus \Y_{(1)}\,,\ 
i_{\V_{(1)}} \colon \X_{(1)}\to \X_{(1)}\oplus \V_{(2)}\oplus \Y_{(1)}\,
\]
are the injections defined by the direct sum decomposition of $\V$ and 
\[
p_{\X_{(1)}} \colon \X_{(1)}\to \X_{(1)}\oplus \V_{(2)}\oplus \Y_{(1)}\,,\ 
p_{\Y_{(1)}} \colon \X_{(1)}\to \X_{(1)}\oplus \V_{(2)}\oplus \Y_{(1)}\,,\ 
p_{\V_{(1)}} \colon \X_{(1)}\to \X_{(1)}\oplus \V_{(2)}\oplus \Y_{(1)}\,
\]
are the corresponding projections.
Set
\begin{multline}\label{Zariski}
\g''=\{z\in \g;\ (z-1)\,,\ A(z)\,,\ C(z)\,,\ (z_h-1)\,,\ A(z_h)\,,\ C(z_h)\,,\ \text{are invertible}\,,\\ E(z)\ne 0\,,\ \text{and}\ E(z_h)\ne 0\}\,.
\end{multline}
This is non-empty Zariski open subset of $\g$. 
From now on we shall assume that $z\in U=\supp\phi\subseteq  \g''$ and let $m_0\in \Bbb Z$ be the minimum of the $m\in \Bbb Z$ such that \eqref{Theorem 7.7.1Hormander4} holds.  We shall impose some more conditions on $\g''$ below.

By using the explicit description of the symplectic form,  \eqref{explicitsymplecticform}, and remembering that the Lie algebra $\g$ acts on $\W$ via minus the right multiplication, we can view the $A=A(z)$, $B=B(z)$, ..., $F=F(z)$ as elements of $\End(\V)$, so that (up to a fixed positive constant multiple relating to $\tr_{\Bbb D/\F}$ and the symplectic form $\langle \cdot,\cdot\rangle$) 
\begin{eqnarray}\label{thephasefunction}
-\phi_{x,x',w_2}(z)&=&\tr_{\Bbb D/\F}\Big( (x-x')^*(x-x')B\\
&+&(x+x'+(x-x')A+w_2E)^*(x+x'+(x-x')A+w_2E)F\nn\\
&+&2(x-x')^*w_2D+w_2^*w_2z_2\Big)\nn\,.
\end{eqnarray}
The derivative (i.e. the linear part) of $-\phi_{x,x',w_2}(z)$, at $z$, viewed as a function of the variables $A$, $B$, $F$, $D$, $E$, ${z_2}$ is given by
\begin{eqnarray*}
\lefteqn{- \phi'_{x,x',w_2}(z)(\Delta_A, \Delta_B, \Delta_F, \Delta_D, \Delta_E, \Delta_{z_2})} \\
&=&\tr_{\Bbb D/\R}\Big((x-x')^*(x-x')\Delta_B\\
&+&((x-x')\Delta_A)^*(x+x'+(x-x')A+w_2E)F\\
&+&(x+x'+(x-x')A+w_2E)^*(x-x')\Delta_AF\\
&+&(w_2\Delta_E)^*(x+x'+(x-x')A+w_2E)F+(x+x'+(x-x')A+w_2E)^*w_2\Delta_EF\\
&+&(x+x'+(x-x')A+w_2E)^*(x+x'+(x-x')A+w_2E)\Delta_F\\
&+&2(x-x')^*w_2\Delta_D+w_2^*w_2\Delta_{z_2}\Big)\,.
\end{eqnarray*}
Notice that $\Delta_AF=F(\Ad(F^{-1})\Delta_A)$.  Also, by the structure of the Lie algebra $\g$, the variables $\Delta_A, \Delta_B, \Delta_F, \Delta_D, \Delta_E, \Delta_{z_2}$ are independent and fill out the corresponding vector spaces. 
The norm $|\phi'_{x,x',w_2}(z)|$ of the functional $\phi'_{x,x',w_2}(z)$, see \cite[Corollary 3, chapter II, paragraph 1]{Weil_Basic} can be estimated from below by taking $\Delta_E=0$ and $\Delta_F=0$. Furthermore, all norms on a finite dimensional vector space are equivalent. 
Hence, with the appropriate choice of the norm $|\cdot |$ on $\End_\Bbb D(\V)$,
\begin{eqnarray}\label{maintechicallemma2}
|\phi'_{x,x',w_2}(z)|
&\geq&|(x-x')^*(x-x')|\\
&+&|(x-x')^*(x+x'+(x-x')A+w_2E)F|\nn\\
&+&|(x+x'+(x-x')A+w_2E)^*(x-x')F\Ad(F^{-1})|\nn\\
&+&2|(x-x')^*w_2|+|w_2^*w_2|\,.\nn
\end{eqnarray}
Using the inequality $|ab|\geq |a||b^{-1}|^{-1}$, which follows from \eqref{operatornorminequality}, and the fact that $|a^*|=|a|$ we see 
that
\begin{eqnarray*}
&&|(x-x')^*(x-x')|\geq |(x-x')^*(x-x')A||A|^{-1}\,,\\
&&|(x-x')^*(x+x'+(x-x')A+w_2E)F|\geq 
|(x-x')^*(x+x'+(x-x')A+w_2E)||F^{-1}|^{-1}\,,\\
&&|(x+x'+(x-x')A+w_2E)^*(x-x')F\Ad(F^{-1})|\\
&&\qquad \geq 
|(x-x')^*(x+x'+(x-x')A+w_2E)|\Ad(F)F^{-1}|^{-1}\,,\\
&&|(x-x')^*w_2|\geq |(x-x')^*w_2E||E|^{-1}\,.
\end{eqnarray*}
Hence,
\begin{eqnarray*}
&&|\phi'_{x,x',w_2}(z)|
\geq C_0(z)\big(|(x-x')^*(x-x')|+|(x-x')^*(x-x')A|+\\
&&\qquad 
|(x-x')^*(x+x'+(x-x')A+w_2E)|+|(x-x')^*w_2E|+|(x-x')^*w_2|+|w_2^*w_2|\big)\,,
\end{eqnarray*}
where
\[
C_0(z)=\min(\frac{1}{2}, C_{00}(z))\,, \ \ \ C_{00}(z)=\min(\frac{1}{2}|A|^{-1}, |C|^{-1}+|\Ad(F)C|^{-1}, |E|^{-1})\,.
\]
Using the triangle 
inequality $|a|+|b|\geq |a\pm b|$ we see that
\begin{eqnarray*}
&&|(x-x')^*(x-x')A|+|(x-x')^*(x+x'+(x-x')A+w_2E)|+
|(x-x')^*w_2E|\\
&\geq&|-(x-x')^*(x-x')A+(x-x')^*(x+x'+(x-x')A+w_2E)-(x-x')^*w_2E|\\
&=&|(x-x')^*(x+x')|.
\end{eqnarray*}
So,
\begin{eqnarray*}
|\phi'_{x,x',w_2}(z)|\geq C_0(z)\big(|(x-x')^*(x-x')|+
|(x-x')^*(x+x')|+|(x-x')^*w_2|+|w_2^*w_2|\big)\,.
\end{eqnarray*}
Our computation applied to $z_h$ shows that
\begin{eqnarray*}
|\phi'_{x,x',w_2}(z_h)|\geq C_0(z_h)(|(x-x')^*(x-x')|+|(x-x')^*(x+x')|+|(x-x')^*w_2|+|w_2^*w_2|)\,.
\end{eqnarray*}
Recall that Lemma \ref{partialweyltransformfor g} provides another expression for the function we would like to estimate, in terms $\phi'_{x,-x',w_2}(z_h)$. Indeed, 
\[
\mathcal K_1(T(\wt{c(z)}))(x,x',w_2)=\det_{\X_1}^{1/2}(\t h)\mathcal K_1(T(\wt{c(z_h)}))(x,-x',w_2)
\]
and as we have seen previously,
\begin{eqnarray*}
|\phi'_{x,-x',w_2}(z_h)|\geq C_0(z_h)(|(x+x')^*(x+x')|+|(x+x')^*(x-x')|+|(x+x')^*w_2|+|w_2^*w_2|)\,.
\end{eqnarray*}
By the triangle inequality,
\begin{eqnarray*}
|(x-x')^*(x-x')|+|(x-x')^*(x+x')|\geq |(x-x'{})^*2x|&=&|(x-x'{})^*x|\,,\\
|(x-x')^*(x'-x)|+|(x-x')^*(x+x')|\geq |(x-x'{})^*2x'|&=&|(x-x'{})^*x'|\,,\\
|(x+x')^*(x+x')|+|(x+x')^*(x-x')|&\geq& |(x+x'{})^*x|\,,\\
|(x+x')^*(x+x')|+|(x+x')^*(x-x')|&\geq& |(x+x'{})^*x'|\,,\\
|(x-x'{})^*x|+|(x+x'{})^*x|&\geq&|x^*x|\,,\\ 
|(-x+x'{})^*x|+|(x+x'{})^*x|&\geq&|x'{}^*x'|\,,\\
|(x-x'{})^*x'|+|(x+x'{})^*x'|&\geq&|x^*x'|\,,\\ |
(-x+x'{})^*x|+|(x+x'{})^*x|&\geq&|x'{}^*x|\,.
\end{eqnarray*}
Therefore
\begin{multline*}
|(x-x')^*(x-x')|+|(x-x')^*(x+x')|+|(x+x')^*(x+x')|+|(x-x')^*(x+x')|\\
\geq \max\{|x^*x|,|x'{}^*x'|,|x'{}^*x|,|x^*x'{}^*|\}\geq \frac{1}{4}
(|x^*x|+|x'{}^*x'|+|x'{}^*x|+|x^*x'{}^*|)\,.
\end{multline*}
Furthermore,
\[
|(x-x')^*w_2|+|(x+x')^*w_2|\geq |x^*w_2|\,,\ \ |(x'-x)^*w_2|+|(x+x')^*w_2|\geq |x'{}^*w_2|\,.
\]
Hence, 
\begin{equation}\label{maintechnicallema00}
\min_{z\in U}|\phi'_{x,x',w_2}(z)|+\min_{z\in U}|\phi'_{x,-x',w_2}(z_h)|
\geq\min(\min_{z\in U}C_0(z), \min_{z\in U}C_0(z_h))\frac{1}{4}m(x,x',w_2)\,,
\end{equation}
where
\[
m(x,x',w_2)=|x^*x|+|x'{}^*x'|+|x^*x'|+|x'{}^*x|+|x^*w_2|+|x'{}^*w_2| +|w_2^*w_2|\,.
\]
Therefore,
\[
\max(\min_{z\in U}|\phi'_{x,x',w_2}(z)|,\min_{z\in U}|\phi'_{x,-x',w_2}(z_h)|)
\geq\frac{1}{8}\min(\min_{z\in U}C_0(z), \min_{z\in U}C_0(z_h))m(x,x',w_2)\,.
\]
We shrink $\g''$ by imposing the additional condition that both $C_0(z)$ and $C_0(z_h)$ are finite and let 
\begin{equation}\label{G''}
\G''=c(\g'')\subseteq \G\,. 
\end{equation}
With the notation of \eqref{thephasefunction} let $f(z)=-\phi_{x,x',w_2}(z)$. Let $x_-=x-x'$, $x_+=x+x'$. 
Then a straightforward computation shows that for $z_0, z, z_0+z\in\g''$
\[
f(z_0+z)=f(z_0)+f'(z_0)(z)+R(z_0,z)(z)(z)\,,
\]
where
\begin{eqnarray*}
R(z_0,z)(z)(z)&=&\tr_{\Df/\F}\big(A^*x_-^*(x_++x_-A_0+w_2E_0)F+E^*w_2^*(x_++x_-A_0+w_2E_0)F\\
&+&(x_++x_-A_0+w_2E_0)^*x_-AF+(x_++x_-A_0+w_2E_0)^*w_2EF\\
&+&A^*x_-^*x_-AF_0+A^*x_-^*w_2EF_0+E^*w_2^*x_-AF_0+E^*w_2^*w_2EF_0\\
&+&A^*x_-^*x_-AF+A^*x_-^*w_2EF+E^*w_2^*x_-AF+E^*w_2^*w_2EF\big)\,,
\end{eqnarray*}
where the subscript $0$ indicates that the corresponding element comes from $z_0$. In order to view this function as in Lemma \ref{Theorem 7.7.1Hormander} we set
\begin{eqnarray*}
R(z_0,z)(y)(y)&=&\tr_{\Df/\F}\big(A_y^*x_-^*(x_++x_-A_0+w_2E_0)F_y+E_y^*w_2^*(x_++x_-A_0+w_2E_0)F_y\\
&+&(x_++x_-A_0+w_2E_0)^*x_-A_yF_y+(x_++x_-A_0+w_2E_0)^*w_2E_yF_y\\
&+&A_y^*x_-^*x_-A_yF_0+A_y^*x_-^*w_2E_yF_0+E_y^*w_2^*x_-A_yF_0+E_y^*w_2^*w_2E_yF_0\\
&+&A^*x_-^*x_-A_yF_y+A^*x_-^*w_2E_yF_y+E^*w_2^*x_-A_yF_y+E^*w_2^*w_2E_yF_y\big)\,,
\end{eqnarray*}
where the subscript $y$ indicates that the corresponding element comes from $y$. 
Since $x_-^*x_+=x^*x+x^*x'-x'{}^*x-x'{}^*x'$, $x_-^*x_-=...$,  it is clear that there is a constant $k(z_0,z)$ depending continuously on $(z_0,z)$ such that
\[
\max_{|y|=1}|R(z_0,z)(y)(y)|\leq k(z_0,z)m(x,x',w_2)\,,
\]
where $m(x,x',w_2)$ is the function defined in \eqref{maintechnicallema00}.
Hence there is a constant $C_U$ such that
\begin{equation}\label{maintechnicallema01}
C_U m(x,x',w_2)\geq\max_{z_0, z_0+z\in U}\max_{|y|=1}|R(z_0,z)(y)(y)|\,.
\end{equation}
Similar analysis applies to $f(z)=\phi_{x,-x',w_2}(z_h)$ with the resulting function $R_h(z_0z)$ so that with the appropriately adjusted constant $C_U$
\[
C_U m(x,x',w_2)\geq\max_{z_0, z_0+z\in U}\max_{|y|=1}|R_h(z_0,z)(y)(y)|\,.
\]
Let $C_\Psi$ be large enough to that for $m(x,x',w_2)>C_\Psi$ 
\[
m(x,x',w_2)> \left(q C_U m(x,x',w_2)\right)^{\frac{1}{2}}
\]
and 
\[
m(x,x',w_2)>q^{m_0}\,,
\]
where $m_0$ is given in Lemma 5.
Then
\begin{multline*}
\max(\min_{z\in U}|\phi'_{x,x',w_2}(z)|,\min_{z\in U}|\phi'_{x,-x',w_2}(z_h)|)\\
\geq 
\max\left( \max_{z_0, z_0+z\in U}\max_{|y|=1}|R(z_0,z)(y)(y)|, \max_{z_0, z_0+z\in U}\max_{|y|=1}|R_h(z_0,z)(y)(y)|, q^{m_0}\right)\,.
\end{multline*}
Therefore
\[
\min_{z\in U}|\phi'_{x,x',w_2}(z)|
\geq 
\max\left( \max_{z_0, z_0+z\in U}\max_{|y|=1}|R(z_0,z)(y)(y)|, q^{m_0}\right)
\]
or
\[
\min_{z\in U}|\phi'_{x,x',w_2}(z_h)|
\geq 
\max\left( \max_{z_0, z_0+z\in U}\max_{|y|=1}|R_h(z_0,z)(y)(y)|, q^{m_0}\right)\,.
\]
Let $\phi(z)=\Psi(\t c(z))\Theta(\wt{c(z)})\gamma(q_{z,\Y_1})$.
By Lemma \ref{Theorem 7.7.1Hormander}, the first condition implies that 
\[
\int_{U}\phi(z)\chi(\phi_{x,x',w_2}(z))\,dz=0
\]
and the second condition that
\[
\int_{U}\phi(z)\chi(\phi_{x,-x',w_2}(z_h))\,dz=0.
\]
But Lemma \ref{partialweyltransformfor g} shows that one expression is a non-zero multiple of the other and the first one is equal to $\mathcal K_1(T(\Psi))(x,x',w_2)$. Hence, \eqref{maintechicallemma1} follows.
\end{prf}
As an immediate consequence of Corollary \ref{Maction} and Proposition \ref{N1action} we deduce the following lemma.
\begin{lem}\label{G'actionandZ'action}
Let $\Zg\subseteq\G$ be the subgroup that acts trivially on $\Y_1^\perp$. Then for $\t n\in \wt{\Zg}$, $v\in \Ss(\X_1,\Ss(\X_2))$, $x_1\in \X_1$ and $\t g'\in\wt{\G'}$,
\begin{equation}\label{G'actionandZ'action1}
\omega(\t n)v(x_1)=\pm \chi_{c(-n)}(2x_1)v(x_1) \,,
\end{equation}
and
\begin{equation}\label{G'actionandZ'action2}
\omega(\t g')v(x_1)=\det_{X_1}^{-1/2}(\t g')\omega_2(\t g')v(g'{}^{-1}x_1) \,.
\end{equation}
\end{lem} 
\section{\bf The functions  $\Psi\in C_c^\infty(\wt\G'')$ act on $\mathcal H_\Pi$ via integral kernel operators. \rm}\label{via the integral kernel operators}\noindent
Given the polarization $\W_2=\X_2\oplus \Y_2$ we have the map
\[
\rho_2:\Ss^*(\W_2)\to \Hom_{\C}(\Ss(\X_2), \Ss^*(\X_2))
\] 
as in \eqref{eqrho}. Then 
\[
1\otimes \rho_2 : \Ss^*(\X_1\times\X_1\times \W_2)\to \Ss^*(\X_1\times\X_1)\otimes \Hom_{\C}(\Ss(\X_2), \Ss^*(\X_2))\,.
\]
In order to shorten the notation we shall write  $\rho_2$ for $1\otimes \rho_2$ and 
\[
\mathcal K_1(T(\t g))(x,x')=\mathcal K_1(T(\t g))(x,x', \cdot)\qquad (x,x'\in \X_1)\,.
\]
In these terms
\begin{equation}\label{complicatedomega}
\omega(\t g)v(x)=\int_{\X_1}\rho_2(\mathcal K_1(T(\t g))(x,x'))(v(x'))\,d\overset . x' \quad (\t g\in\wt\G, v\in \Ss(\X_1,\Ss(\X_2)); x,x'\in \X_1)\,,
\end{equation}
i.e. $\omega(\t g)v$ is a function taking $x\in \X_1$ to $\Ss(\X_2)$.
Let $\X_1^{max}\subseteq \X_1 = \Hom_{\Bbb D}(\X_{(1)}, \V')$ be the subset of the surjective maps. The stable range assumption implies that this is a dense subset. 
Let $\Psi\in C_c^\infty(\wt\G'')$ as in Lemma \ref{maintechicallemma}.
For fixed $x, x'\in \X_1^{max}$ the operator norm of 
\begin{equation}\label{hsnorm}
\rho_2(\mathcal K_1(T(\Psi))(x,x') \ \in \Hom_{\C}(\L^2(\X_2), \L^2(\X_2))
\end{equation}
is bounded by the Hilbert-Schmidt norm, which is finite. Indeed, Lemma \ref{maintechicallemma} shows that
\[
\mathcal K_1(T(\Psi))(x,x', w_2)
\]
is a compactly supported  function of $x^*w_2$ and hence of $w_2$, because $x^*$ is an injective map from $\V'$ to $\X_{1}$. Therefore
\[
\mathcal K_1(T(\Psi))(x,x', \cdot)\in \L^2(\W_2)\,,
\]
which means that the Hilbert-Schmidt norm of \eqref{hsnorm} is finite.

In general, we denote by $\sigma^c$ the representation contragredient to $\sigma$ and by $\mathcal H_\sigma$  a Hilbert space where $\sigma$ is realized.

The group $\G'$ acts on $\X_1^{max}$, via the left multiplication, so that the quotient $\G'\backslash \X_1^{max}$ is a manifold. If $dx$ is a Lebesgue measure on $\X_1$, we shall denote by $d\overset . x$ the corresponding quotient measure on $\G'\backslash \X_1^{max}$.
Let $\mathcal U$ be the Hilbert space of functions $u:\X_1^{max}\to \L^2(\X_2)\otimes \mathcal H_{\Pi'{}^c}$ such that for all $\t g'\in \wt\G'$
\begin{eqnarray}\label{realization of Pi}
u(g'x)=(\omega_2\otimes \det_{\X_1}^{-1/2}\Pi'{}^c)(\t g')u(x)\ \ \ \text{and}\ \ \ 
\int_{\G'\backslash \X_1^{max}}\parallel u(x)\parallel^2\,d\overset . x<\infty\,,
\end{eqnarray}
where $\det_{\X_1}^{-1/2}$ is as in \eqref{G'actionandZ'action2}.
\begin{lem}\label{relizationasintegralkernels}
The representation $\Pi$ is realized on the Hilbert space $\mathcal U$ define in \eqref{realization of Pi}.
For $\Psi\in C_c^\infty(\wt\G'')$, the operator $\Pi(\Psi)$ is given in terms of an integral kernel defined on $\X_1^{max}\times\X_1^{max}$ as follows
\[
(\Pi(\Psi)u)(x)=\int_{\G'\backslash \X_1^{max}}K_\Pi(\Psi)(x,x')u(x')\,d\overset . x' \quad(u\in \mathcal H_\Pi)\,,
\]
where
\begin{equation}\label{relizationasintegralkernels1}
K_\Pi(\Psi)(x,x')=\int_{\G'}\omega_2(\t g) \rho_2 (\mathcal K_1 (T(\Psi))(g^{-1}x,x',\cdot))\otimes \det_{\X_1}^{-1/2}(\t g)\Pi'{}^c(\t g)\,dg\,.
\end{equation}
Furthermore,
\begin{eqnarray}\label{relizationasintegralkernels2}
\lefteqn{\tr K_\Pi(\Psi)(x,x')} \\
&=&\int_{\G'} \int_{\W_2} T_2(\t g)(w_2) \mathcal K_1 (T(\Psi))(g^{-1}x,x', w_2)\det_{\X_1}^{-1/2}(\t g)\Theta_{\Pi'{}^c}(\t g)\,dw_2\,dg\,,\nn
\end{eqnarray}
where $\int_{\W_2} T_2(\t g)(w_2)\phi(w_2)\,dw_2$ stands for $T_2(\t g)(\phi)$. 
\end{lem}
\begin{prf}
We proceed as in \cite[Proposition 4.8]{DaszPrzebindaInv}. Define a map 
\[
\Qg: \Ss(\X_1, \Ss(\X_2))\otimes \mathcal H_{\Pi'{}^c}\to {\mathcal U}
\] 
by
\begin{eqnarray}\label{realization of Q}
\Qg(v\otimes \eta)(x)=\int_{\G'}(\omega\otimes\Pi'{}^c)(\t g)(v (x) \otimes \eta)\, dg\,  \quad  (x \in \X_1^{max}).
\end{eqnarray}
Then \eqref{G'actionandZ'action2} shows that
\[
\Qg(v\otimes \eta)(x)=\int_{\G'}\omega_2(\t g)(v(g^{-1}x))\otimes\det_{\X_1}^{-1/2}(\t g)\Pi'{}^c(\t g)\eta\, dg\,.
\]
This last integral converges because $|g^{-1}x|$ is a constant multiple of the norm of $g$, as defined in \cite[2.A.2.4]{WallachI}. (The constant depends on $x$, which is fixed.) 
The argument used in the proof of Lemma 3.11 in \cite{DaszPrzebindaInv} shows that the range of $\Qg$ is dense in $\mathcal U$. 
The action of  $\t g\in \wt\G$ on $\mathcal U$ is defined via the the action of $\omega(\t g)$ on the $v$.  We denote $\pi(\t g)=\det_{\X_1}^{-1/2}(\t g)\Pi'{}^c(\t g)$ and have
\begin{eqnarray*}
\lefteqn{\Qg(\omega(\Psi)v\otimes\eta)(x)} \\
&=&\int_{\G'}\omega_2(\t g)((\omega(\Psi)v)(g^{-1}x))\otimes \pi(\t g)\eta\,d g\\
&=&\int_{\G'}\int_{\X_1^{max}}\omega_2(\t g)\rho_2(\mathcal K_1(T(\Psi))(g^{-1}x,x'))(v(x'))\otimes \pi(\t g)\eta\,dx'\,d g  \quad \text (By\  \eqref{complicatedomega}.) \\
&=&\int_{\X_1^{max}}\int_{\G'}\omega_2(\t g)\rho_2(\mathcal K_1(T(\Psi))(g^{-1}x,x'))(v(x'))\otimes \pi(\t g)\eta\,d g\,dx'\\
&=&\int_{ \G'\backslash\X_1^{max}}\int_{\G'}\int_{\G'}\omega_2(\t g)\rho_2(\mathcal K_1(T(\Psi))(g^{-1}x,h^{-1}x'))(v(h^{-1}x'))\otimes \pi(\t g)\eta\,d g\,dh\,d\overset . x'\\
&=&\int_{ \G'\backslash\X_1^{max}}\int_{\G'}\int_{\G'}\omega_2(\t g\t h)\rho_2(\mathcal K_1(T(\Psi))(h^{-1}g^{-1}x,h^{-1}x'))(v(h^{-1}x'))\otimes \pi(\t g\t h)\eta\,d g\,dh\,d\overset . x'\\
&=&\int_{ \G'\backslash\X_1^{max}}\int_{\G'}\int_{\G'}\omega_2(\t g\t h)\omega_2(\t h)^{-1}\rho_2(\mathcal K_1(T(\Psi))(g^{-1}x,x'))(\omega_2(\t h)v(h^{-1}x'))\otimes \pi(\t g\t h)\eta\,d g\,dh\,d\overset . x'\\
&=&\int_{ \G'\backslash\X_1^{max}}\left(\int_{\G'}\omega_2(\t g)\rho_2(\mathcal K_1(T(\Psi))(g^{-1}x,x'))\otimes \pi(\t g)\,dg\right)\Qg(v\otimes\eta)(x')\,d\overset . x'\,,
\end{eqnarray*}
where by Lemma \ref{maintechicallemma} all the integrals are convergent.
This verifies \eqref{relizationasintegralkernels1}. 

Furthermore, the usual argument shows that $\mathcal K_1(T(\Psi))(g^{-1}x,x',w_2)$ is a differentiable function of $g$ and $w_2$ so that the method of stationary phase applies to ensure the convergence of the integrals below,
\begin{eqnarray}\label{lastequality}
\tr K_\Pi(\Psi)(x,x')
&=&\tr \int_{\G'}\omega_2(\t g) \rho_2 (\mathcal K_1 (T(\Psi))(g^{-1}x,x',\cdot)\otimes \det_{\X_1}^{-1/2}(\t g)\Pi'{}^c(\t g)\,dg\\
&=&\int_{\G'}\int_{\W_2}T_2(\t g)(w_2) (\mathcal K_1 (T(\Psi))(g^{-1}x,x',w_2)\det_{\X_1}^{-1/2}(\t g)\Theta_{\Pi'{}^c}(\t g)\,dg\,d w_2\,.\nn
\end{eqnarray}
Therefore \eqref{relizationasintegralkernels2} follows from \eqref{relizationasintegralkernels1}. 

Here is an explanation of the second equality in \eqref{lastequality}.
If $\W_2=0$ then $\mathcal K_1=\mathcal K$, \\ $\mathcal K_1 (T(\Psi))(g^{-1}x,x',\cdot)=\mathcal K_1 (T(\Psi))(g^{-1}x,x')$ is a scalar $C_c^\infty$ function of $g$. Therefore, by Harish-Chandra
\begin{multline*}
\tr K_\Pi(\Psi)(x,x')
=\tr \int_{\G'} \mathcal K (T(\Psi))(g^{-1}x,x') \det_{\X}^{-1/2}(\t g)\Pi'{}^c(\t g)\,dg\\
=\int_{\G'} \mathcal K (T(\Psi))(g^{-1}x,x')\det_{\X}^{-1/2}(\t g)\Theta_{\Pi'{}^c}(\t g)\,dg\,.
\end{multline*}
Suppose $\W_2\ne 0$. 
Then
\begin{multline*}
\tr \int_{\G'}\omega_2(\t g) \rho_2 (\mathcal K_1 (T(\Psi))(g^{-1}x,x',\cdot)\otimes \det_{\X_1}^{-1/2}(\t g)\Pi'{}^c(\t g)\,dg\\
=\tr \int_{\G'}\tr(\omega_2(\t g) \rho_2 (\mathcal K_1 (T(\Psi))(g^{-1}x,x',\cdot))\det_{\X_1}^{-1/2}(\t g)\Pi'{}^c(\t g)\,dg\\
=\tr \int_{\G'}\int_{\W_2}T_2(\t g)(w_2) \mathcal K_1 (T(\Psi))(g^{-1}x,x',w_2)\,dw_2\det_{\X_1}^{-1/2}(\t g)\Pi'{}^c(\t g)\,dg\,,
\end{multline*}
where the last equality follows from \eqref{eqtromegarho}. This coincides with the last expression in \eqref{lastequality}.
\end{prf}
\section{\bf The equality $\Theta_\Pi=\Theta'_{\Pi'}$  on a non-empty Zariski open subset of $\wt \G''\subseteq \wt \G$. \rm}\label{The character}\noindent
For $\phi, \psi \in \Ss(\W_2)$, we set
\[
\phi \natural \psi(w') = \int_{\W_2} \phi(w) \psi(w'-w) \chi(\dfrac{1}{2} \langle w, w' \rangle) dw. 
\]
Assume that the union of the conjugacy classes of Cartan subgroups in $\G'$ is dense in $\G'$. 
Here the subset $\G''\subseteq \G$ was defined in \eqref{G''}. Fix a test function $\Psi\in \C_c^\infty(\wt\G'')$. Lemma \ref{maintechicallemma} implies that all the consecutive integrals in the following computation are absolutely convergent:
\begin{eqnarray}\label{firstcharacterformula1}
\lefteqn{\Theta_\Pi(\Psi)=\tr \Pi(\Psi)
=\int_{\G'\backslash\X_1^{max}}\tr { K_\Pi(\Psi)(x,x)} \,d\overset . x \quad \text{ (By Lemma \ref{relizationasintegralkernels}.)} } \nn \\
&=&\int_{\G'\backslash\X_1^{max}}\tr { K_\Pi(\Psi)(-x,-x) }\,d\overset . x\nn\\
&=&\int_{\G'\backslash\X_1^{max}}\int_{\G'}\int_{\W_2}T_2(\t g)(w_2) \mathcal K_1(T(\Psi))(-g^{-1}x,-x,w_2) \det_{\X_1}^{-1/2}(\t g)\Theta_{{\nn\Pi'}}(\t g^{-1})\,dw_2\,dg\, d\overset . x\nn\\
&=&\chi_{{\Pi'}}((-1\t ))\int_{\G'\backslash\X_1^{max}}\int_{\G'}\int_{\W_2}T_2((-1\t )\t g)(w_2) \mathcal K_1(T(\Psi))(g^{-1}x,-x,w_2)\nn\\
&&\det_{\X_1}^{-1/2}((-1\t )\t g)\Theta_{{\Pi'}}(\t g^{-1})\,dw_2\,dg\, d\overset . x\nn\\
&=&\chi_{{\Pi'}}((-1\t ))\int_{\G'\backslash\X_1^{max}}\int_{\G'}\left(T_2((-1\t ))\natural T_2(\t g)\natural \mathcal K_1(T(\Psi))(g^{-1}x,-x,\cdot)\right)(0)\nn\\
&&\det_{\X_1}^{-1/2}((-1\t )\t g)\Theta_{{\Pi'}}(\t g^{-1})\,dg\, d\overset . x\nn\\
&=&\chi_{{\Pi'}}((-1\t ))\Theta_2((-1\t ))\int_{\G'\backslash\X_1^{max}}\int_{\G'}\int_{\W_2}\left(T_2(\t g)\natural \mathcal K_1(T(\Psi))(g^{-1}x,-x,\cdot)\right)(w_2)\nn\\
&&\det_{\X_1}^{-1/2}((-1\t )\t g)\Theta_{{\Pi'}}(\t g^{-1})\,dw_2\,dg\, d\overset . x\nn\\
&=&\chi_{{\Pi'}}((-1\t ))\Theta_2((-1\t ))\det_{\X_1}^{-1/2}((-1\t ))\int_{\G'\backslash\X_1^{max}}\int_{\G'}\int_{\W_2}\mathcal K_1(T_2(\t g)\natural T(\Psi))(g^{-1}x,-x, w_2)\nn\\
&&\det_{\X_1}^{-1/2}(\t g)\Theta_{{\Pi'}}(\t g^{-1})\,dw_2\,dg\, d\overset . x\,.
\end{eqnarray}
Recall that with the appropriate notion of the tensor product we have $\Ss(\W)=\Ss(\W_1)\otimes\Ss(\W_2)$. It is easy to check that for $\phi_1, \psi_1\in \Ss(\W_1)$ and $\psi_2\in \Ss(\W_2)$,
\[
(\phi_1\otimes \delta)\natural (\psi_1\otimes \psi_2)=(\phi_1\natural\psi_1)\otimes\psi_2\,,
\]
where the $\natural$ on the right hand side happens in $\Ss(\W_1)$. Hence, for $x,x'\in\X_1$ and $w_2\in\W_2$,
\begin{multline*}
\mathcal K_1((\phi_1\otimes \delta)\natural (\psi_1\otimes \psi_2))(x,x',w_2)
=\int_{\X_1}\mathcal K_1(\phi_1)(x,x'')\mathcal K_1(\psi_1)(x'',x')\,dx''\psi(w_2)\\
=\int_{\X_1}\mathcal K_1(\phi_1)(x,x'')\mathcal K_1(\psi_1\otimes\psi_2)(x'',x',w_2)\,dx''\,.
\end{multline*}
Hence by a linear approximation, for any $\psi\in\Ss(\W)$, 
\[
\mathcal K_1((\phi_1\otimes \delta)\natural \psi)(x,x',w_2)
=\int_{\X_1}\mathcal K_1(\phi_1)(x,x'')\mathcal K_1(\psi)(x'',x',w_2)\,dx''\,.
\]
For $g \in \GL(\X)$ we can use $\phi_1$ to approximate $T_1(\t g)$. Since $T_1(\t g)$ is identified with $T_1(\t g)\otimes \delta$, Proposition \ref{formula for M R} shows that
\begin{multline*}
\mathcal K_1((T_1(\t g))\natural \psi)(x,x',w_2)
=\int_{\X_1}\det_{\X_1}^{-1/2}(\t g)\delta(g^{-1}x-x'')\mathcal K_1(\psi)(x'',x',w_2)\,dx''\\
=\det_{\X_1}^{-1/2}(\t g)\mathcal K_1(\psi)(g^{-1}x,x',w_2)\,.
\end{multline*}
Now we substitute $T_2(\t g)\natural T(\Psi)$ for $\psi$ and and $-x$ for $x'$ to see that
\[
\mathcal K_1(T_1(\t g)\natural T_2(\t g)\natural T(\Psi))(x,-x, w_2) 
= \det_{\X_1}^{-1/2}(\t g)\mathcal K_1(T_2(\t g)\natural T(\Psi))(g^{-1}x,-x, w_2)\,.
\]
Since $T_1(\t g)\natural T_2(\t g)=T(\t g)$,  \eqref{firstcharacterformula1} is equal to
\begin{eqnarray*}
&&\chi_{{\Pi'}}((-1\t ))\Theta_2((-1\t ))\det_{\X_1}^{-1/2}((-1\t ))
\int_{\G'\backslash\X_1^{max}}\int_{\G'}\int_{\W_2}\int_{\Y_1} T(\t g)\natural T(\Psi)(x+x+y+w_2)\nn\\
&&\chi(\frac{1}{2}\langle y, x-x\rangle)\Theta_{{\Pi'}}(\t g^{-1})\,dy\,dw_2\,dg\, d\overset . x.
\end{eqnarray*}
and finally, the change of variable $x$ to $\frac{1}{2}x$ gives
\begin{eqnarray}
&=&\chi_{{\Pi'}}((-1\t ))\Theta((-1\t ))
\int_{\G'\backslash\X_1^{max}}\int_{\G'}\int_{\W_2}\int_{\Y_1} T(\t g)\natural T(\Psi)(x+y+w_2)\nn\\
&&\Theta_{{\Pi'}}(\t g^{-1})\,dy\,dw_2\,dg\, d\overset . x \, , \label{firstcharacterformula1a}
\end{eqnarray}
where the function under the integral is constant on the fibers of the covering map because we assume that $\Pi'$ is genuine. Also, the integral over $(\G'\backslash\X_1^{max})\times \G'$ is also absolutely convergent. 
For a function $f$ supported in $Z'\G'{}^\circ$, we apply the Weyl - Harish-Chandra integration formula for $\G'$
\begin{equation}\label{Weyl - Harish-Chandra integration formula}
\int_{\G'}f(g')\,dg'=\sum_{\H'}\frac{1}{|W(\H')|}\int_{\reg{\H'}}\int_{\G'/\H'}f(g'h'g'{}^{-1})\,d\overset . g'\,|\Delta(h')|^2\,dh'
\end{equation}
(see appendix \ref{appenA})
to the integral over $\G'$ in  
\eqref{firstcharacterformula1a} and see that 
\begin{eqnarray}\label{firstcharacterformula2}
&&\Theta_\Pi(\Psi)\\
&=&\chi_{\Pi'}((-1\t ))\Theta((-1\t ))
\sum_{\H'}\frac{1}{|W(\H')|}\int_{\G'\backslash\X_1^{max}}\int_{\reg{\H'}}\int_{\G'/\H'}\int_{\W_2}\int_{\Y_1} T(g'\t h'g'{}^{-1})\natural T(\Psi)(x+y+w_2)\nn\\
&&\Theta_{\Pi'}(\t h'{}^{-1})\,|\Delta(h')|^2\,dy\,dw_2\,d\overset . g'\,dh'\, d\overset . x\nn\\
&=&\chi_{\Pi'}((-1\t ))\Theta((-1\t ))
\sum_{\H'}\frac{1}{|W(\H')|}\int_{\H'\backslash\X_1^{max}}\int_{\reg{\H'}}\int_{\W_2}\int_{\Y_1} T(\t h')\natural T(\Psi)(x+y+w_2)\nn\\
&&\Theta_{\Pi'}(\t h'{}^{-1})\,|\Delta(h')|^2\,dy\,dw_2\,dh'\, d\overset . x\nn\\
&=&\chi_{\Pi'}((-1\t ))\Theta((-1\t ))
\sum_{\H'}\frac{1}{|W(\H')|}\int_{\H'\backslash\X_1^{max}}\int_{\W_2}\int_{\Y_1}  \int_{\reg{\H'}} T(\t h')\natural T(\Psi)(x+y+w_2)\nn\\
&&\Theta_{\Pi'}(\t h'{}^{-1})\,|\Delta(h')|^2\,dh'\, dy\,dw_2\, d\overset . x\nn\\
&=&\chi_{\Pi'}((-1\t ))\Theta((-1\t ))
\sum_{\H'}\frac{1}{|W(\H')|} \int_{\H'\backslash\W^{max}} \int_{\reg{\H'}}  T(\t h')\natural T(\Psi)(w)\nn\\
&&\Theta_{\Pi'}(\t h'{}^{-1})\,|\Delta(h')|^2\, dh' \, d\overset. w\nn\\
&=&\chi_{\Pi'}((-1\t ))\Theta((-1\t ))
\sum_{\H'}\frac{1}{|W(\H')|}\int_{\reg{\H'}}\Theta_{\Pi'}(\t h'{}^{-1})\,|\Delta(h')|^2\\
&&\int_{\H'\backslash\W^{max}}\int_{\wt\G}\Psi(\t g)T(\t h'\t g)(w)
\,d\t g\, d\overset . w\,dh'\nn\,.
\end{eqnarray}
Let $\A'$ be the $\F$-split component of $\H'$. Fix $h'\in\reg{\H'}$ and notice that
\begin{multline*}
\int_{\A'\backslash\W^{max}}\int_{\wt\G}\Psi(\t g)T(\t h'\t g)(w)
\,d\t g\, d\overset . w
=\int_{\H'\backslash\W^{max}}\int_{\A'\backslash\H'}\int_{\wt\G}\Psi(\t g)T(\t h'\t g)(h_1w)
\,d\t g\,d\overset . h_1\, d\overset . w\\
=\int_{\H'\backslash\W^{max}}\int_{\A'\backslash\H'}\int_{\wt\G}\Psi(\t g)T(h_1^{-1}\t h'h_1\t g)(w)
\,d\t g\,d\overset . h_1\, d\overset . w\\
=\int_{\H'\backslash\W^{max}}\mathrm{volume}(\A'\backslash\H')\int_{\wt\G}\Psi(\t g)T(\t h'\t g)(w)
\,d\t g\, d\overset . w\,.
\end{multline*}
Hence 
\begin{multline*}
\int_{\H'\backslash\W^{max}}\int_{\wt\G}\Psi(\t g)T(\t h'\t g)(w)
\,d\t g\, d\overset . w
=\frac{1}{\mathrm{volume}(\A'\backslash\H')}
\int_{\A'\backslash\W^{max}}\int_{\wt\G}\Psi(\t g)T(\t h'\t g)(w)
\,d\t g\, d\overset . w\,.
\end{multline*}
The integral on the right hand side
\[
\int_{\A'\backslash\W^{max}}\int_{\wt\G}\Psi(\t g)T(\t h'\t g)(w)
\,d\t g\, d\overset . w
=\int_{\A'\backslash\W^{max}}\int_{\wt\G}\Psi(\t h'^{-1}\t g)T(\t g)(w)
\,d\t g\, d\overset . w
\]
 extends by the same formula to $\Psi\in C_c^\infty(\wt {\A''{}^c})$ as in \cite[(127)]{LokePrzebinda_chc_padic_def}. By \cite[Proposition 27]{LokePrzebinda_chc_padic_def} this extension has a unique restriction to $\wt{\G}$. Thus
\[
\int_{\A'\backslash\W^{max}}\int_{\wt\G}\Psi(\t g)T(\t h'\t g)(w)
\,d\t g\, d\overset . w
=\Chc _{\t h'}(\Psi)\,.
\]
Therefore \eqref{firstcharacterformula2} shows that
\begin{multline*}
\Theta_\Pi(\Psi)=\\
\chi_{\Pi'}((-1\t ))\Theta((-1\t ))
\sum_{\H'}\frac{1}{|W(\H')|}\int_{\reg{\H'}}\Theta_{\Pi'}(\t h'{}^{-1})\,|\Delta(h')|^2\frac{1}{\mathrm{volume}(\A'\backslash\H')}\Chc _{\t h'}(\Psi)
\,dh'\\
=\Theta'_{\Pi'}(\Psi)\,,
\end{multline*}
Here the last equality is \eqref{Conjectureforcharacters2}
This completes the proof of Theorem \ref{theteequalstheta'}.


\appendix

\section{\bf The Weyl - Harish-Chandra integration formula}\label{appenA}
\setcounter{thh}{0}
\renewcommand{\thethh}{A.\fontindex{thh}}
\setcounter{equation}{0}
\renewcommand{\theequation}{A.\fontindex{equation}}
In order to have the formula \eqref{Weyl - Harish-Chandra integration formula} we need to know that the union of the conjugacy classes of all Cartan subgroups is dense in the group.
Let $\overline{\F}$ be the algebraic closure of~$\F$.
Let $\E$ be a connected reductive algebraic group defined over $\F$. 
We define a Cartan subgroup of $\E(\F)$ as the centralizer of a Cartan subalgebra of the Lie algebra  of $\E(\F)$.

First we suppose $\E$ is a connected semisimple algebraic group defined over $\F$.
An element $g \in \E(\overline{\F})$ is \emph{strongly regular} if its centralizer in $\E(\overline{\F})$ is a Cartan subgroup.
Since $\F$ is a perfect field, the set $R(\overline{\F})$ of strongly regular elements in~$\E(\overline{\F})$ is a Zariski open dense subset by \cite[Section 2.15]{Steinberg_Regular_elements}. Let $W(\overline{\F}) = \E(\overline{\F}) \setminus R(\overline{\F})$. This is a proper Zariski closed subset of $\E(\overline{\F})$ defined over $\F$. 

\begin{lem}
The intersection $R(\F)=R(\overline{\F})\cap \E(\F)$ is a dense subset of $\E(\F)$ in the $\F$-analytic topology.
\end{lem}

\begin{proof}
By Proposition 2.5.2 in \cite{Margulis_book_Discrete_subgroups}, $\E(\F)$ is a pure $\F$-analytic manifold.
Suppose $R(\F)$ is not dense. Then there is a point $x \in W(\F) = \E(\F) \setminus R(\F)$ and an $\F$-analytic open subset $U_x$ of $\E(\F)$ such that $x \in U_x$ and $U_x \cap R(\F) = \emptyset$. Hence $U_x \subseteq W(\F)$.
Let $\mu$ be a Haar measure on $\E(\F)$. Then $\mu(U_x) > 0$ so $\mu(W(\F)) > 0$. Now $W = \bigcup W_i$ is a finite union of irreducible subvarieties. This contradicts Proposition 2.5.3(i) in \cite{Margulis_book_Discrete_subgroups} which states that $\mu(W_i(F)) = 0$.
\end{proof}

We return to our original setting where $\E$ is a connected reductive algebraic group defined over $\F$.
Let $S$ be the union of conjugates of Cartan subgroups of $\E(\F)$. 

\begin{pro} \label{propdense}
The subset $S$ is a dense subset of $\E(\F)$ in the $\F$-analytic topology.
\end{pro}

\begin{prf}
We write $\E(\F) = C \E^s(\F)$ where $C$ is the center of $\E(\F)$ and $\E^s(\F)$ is a the semisimple subgroup. 
By the last lemma $R(\F)$ is a dense subset of $\E^s(\F)$. Hence $C R(\F)$ is a dense subset of $\E(\F)$. The proposition follows because $S$ contains~$C R(\F)$.
\end{prf}

\subsection{}
A member of an irreducible dual pair in a symplectic group defined over $\F$ is a the set  of the $\F$ rational points $\E(\F)$ of an algebraic group $\E$ defined over $\F$, which is connected except the case of an odd orthogonal group.
A Cartan subgroup of $\E(\F)$ is equal to a Cartan subgroup of the identity component together with the center $Z$ of the group.
Let $\widetilde{\E}(\F)$ be the preimage of $\E(\F)$ in the metaplectic group. 
Then $\widetilde{\E}(\F)$ is a double cover of  $\E(\F)$.
Over the algebraic closure $\overline{\F}$, $\E$ is either a symplectic group, an orthogonal group or a general linear group. 
Then a Cartan subgroup of $\widetilde{\E}(\F)$ is the preimage of a Cartan subgroup  of $\E(\F)$. The union of all of them is equal to the preimage $\widetilde{SZ} \subseteq \widetilde{\E}(\F)$ of $SZ$.

\begin{pro}
The subset $\widetilde{SZ}$ is dense in $\widetilde{\E}(\F)$ if and only if $\E$ is not an even orthogonal group.
\end{pro}

The above proposition is equivalent to the next proposition.

\begin{pro}
The subset $SZ$ is dense in $\E(\F)$ if and only if $\E$ is not an even orthogonal group.
\end{pro}

\begin{prf}
If $\E$ is the symplectic group or the general linear group or an odd orthogonal group, the proposition follows from Proposition~\ref{propdense} by setting $S = SZ$.


If $\E$ is an odd orthogonal group.
Then $\E(\F) = \E^\circ(\F) \times \{ \pm 1 \}$ where $\E^\circ(\F)$ is the special orthogonal group.
Let $S$ be the union of conjugates of Cartan subgroups of $\E^\circ(\F)$.
Then
\[
SZ = S \times \{ \pm 1 \}.
\]
By Proposition \ref{propdense}, $S$ is dense in $\E^\circ(\F)$ so $SZ$ is dense in $\E(\F)$.

If $\E$ is an even orthogonal group. 
Then $\E(\overline{\F})$ is isomorphic to $\Og(2n,\overline{\F})$, the split orthogonal group defined by split symmetric form
\[
x_1 x_{n+1} + x_2 x_{n+2} + \ldots + x_n x_{2n}
\]
on $\overline{\F}^{2n}$. Let $\T(\F)$ be a Cartan subgroup $\E(\F)$. 
We may further assume that $\T(\overline{\F})$ is the subgroup of diagonal matrices. 
Then $\T(\overline{\F}) \subseteq \SO(2n,\overline{\F})$. In particular $S \subseteq \SO(2n,\overline{\F})$ and $S$ cannot be dense in $\E(\F)$.
\end{prf}

\biblio
\end{document}